\documentclass[10pt, journal]{IEEEtran}

\usepackage{amssymb}
\usepackage{amsmath}
\usepackage{multicol}
\usepackage[numbers]{natbib} 
\usepackage{balance}
\usepackage[pdftex]{graphicx}
\usepackage{subcaption}

\usepackage{multirow}
\usepackage{longtable}
\usepackage{tabularx}

\usepackage{xcolor}
\usepackage{url}
\usepackage{arydshln}
\begin{document}
\title{EV Fleet Flexibility Estimation and Forecasting for V2X Applications}


\author{Chaimaa~Essayeh,
        Amin~Vilan,
        Omid~Homaee,      and~Vahid~Vahidinasab,~\IEEEmembership{Senior~Member,~IEEE}
\thanks{C. Essayeh and O. Homaee  are with the Department
of Engineering, Nottingham Trent University, United Kingdom.}
\thanks{A. Vilan is with the Islamic Azad University, Tehran, Iran}
\thanks{V. Vahidinasab is with the Business School, University of Salford}
\thanks{Corresponding author e-mail: v.vahidinasab@salford.ac.uk}
\thanks{Manuscript received Month xx, yyyy; revised Month xx, yyyy.}}

\markboth{Journal of \LaTeX\ Class Files,~Vol.~14, No.~8, August~2015}%
{Shell \MakeLowercase{\textit{et al.}}: Bare Demo of IEEEtran.cls for IEEE Journals}

\maketitle 

\begin{abstract}
Forecasting the flexibility potential of Vehicle-to-Everything (V2X) systems is important for the future of energy networks, where the integration of renewable energy sources and electric vehicles poses significant challenges. In this paper, we present a novel method for estimating and predicting V2X flexibility potential of an EV fleet, based on an aggregate polytope representation, addressing the need for accurate and reliable forecasting methods in the realm of sustainable transportation. The method is robust against individual uncertainties of EV owners behaviours as it is applied at an aggregate level, and  the reformulation of the V2X potential as a set of linear constraints allows the proposed method to be integrated into different optimisation problems and therefore be applied for diverse V2X applications. Case studies showcase the capability of the method in capturing the V2X flexibility potential and demonstrate its effectiveness for different V2X applications.
\end{abstract}


\IEEEpeerreviewmaketitle


\section{Introduction}
\label{introduction}
The transportation sector, one of the largest contributors to greenhouse gas emissions \cite{zhang2020role} due to its major reliance on fossil fuels, is undergoing a significant transition driven by government policies aimed at reducing climate change and pollution \cite{PEREIRINHA2018235}. Electrification of transportation, particularly through Battery Electric Vehicles (BEVs) powered by renewable energy, is seen as a key solution for cutting emissions and improving air quality in cities \cite{FERNANDEZ2019476}. While the shift is still in its early stages, advancements in battery technology have resulted in larger capacities and lower costs, making EVs more accessible. This trend is expected to accelerate, with countries like the UK projected to have up to 28 million EVs on the road by 2035 \cite{EnergySavingTrust}.

 This level of electrification will significantly impact the operation of electrical energy systems. The primary challenge will be the substantial increase in electrical energy demand required to supply EVs. For instance, the transportation sector accounts for roughly 36\% of energy consumption in the UK. Transitioning this energy demand to electrical networks would create considerable challenges related to generation, transmission, and distribution \cite{5409638}. From a more technical point of view, this could lead to voltage stability and power quality issues, which have been thoroughly reviewed in  \cite{IET_EV, HAIDAR2014689}, and may necessitate the deployment of costly grid reinforcement.

Although this transition in the transportation sector presents significant challenges for the energy sector, it also offers several opportunities. EV batteries can serve as mobile energy storage systems that can be used to facilitate the integration of intermittent renewable energy resources \cite{LUND20083578, LIU2013445}, support the system in high demand load periods \cite{YUAN2021121564} and supply backup power during outages. Unlocking the EV flexibility potential offers a cost-effective solution to grid reinforcement, generates new revenue streams for end-customers and enhances EVs'role in modernizing and stabilizing the energy system.
To effectively harness the benefit of this flexibility, it is essential to develop various flexibility  mechanisms. Particularly, and due to the uncertain behaviour of EVs' owners, accurate EV flexibility estimation and forecasting methods will play a pivotal role in the adoption of such flexibility. 
%


To provide a reliable foundation for flexibility forecasting, robust estimation methods are essential. In \cite{flex_south_korea}, flexibility is estimated by multiplying the charging station's utilization rate, the power capacity of the charging ports, and the service duration. This method assumes that the operator can draw power equivalent to the charging port's capacity from each connected EV. However, if an EV has just been plugged in with a minimal state of charge ($SOC_{min}$), extracting V2X power could push it below the recommended minimum level. Likewise, if the EV has already reached its maximum state of charge ($SOC_{max}$), the approximate formula implies that the EV can be charged beyond this limit, simply because it is plugged in and theoretically capable of being charged or discharged at the port's power rate. An enhanced estimation method was used in \cite{7778827}, that instead of using the service duration to estimate the flexibility used the idle period. The study estimated flexibility by computing the required charging time to reach the desired energy level and subtracting it from the EV's availability period, assuming the remaining idle period can be used for flexibility. 
However, this assumption overlooks user preferences and technical constraints. For example, if an EV reaches full charge at time $t$ and is scheduled to depart at $t+1$, the method incorrectly considers this final time slot as available for flexibility. Charging is not possible because the battery is already full, and discharging would leave the EV with less than the desired charge level, conflicting with the user's needs. Thus, this approach risks overestimating the actual flexibility potential.
%
%
%
In their approach, \cite{BARTOLINI2023101196} presented a more accurate model for flexibility estimation. The flexibility potential of a carpark is estimated by first calculating the baseline schedule for each EV and then using linear programming to compute the available flexibility over a specific time window, followed by summing the individual flexibilities to estimate the aggregate flexibility. However, this method presents several limitations. Assuming the same baseline for all EVs is unrealistic, as individual behaviors vary significantly. Moreover, computing flexibility for each EV in every time window requires re-running the procedure for all windows, which becomes computationally expensive. 
%
%
The studies in \cite{7463483, 9837819}  modeled the power and energy boundaries of an entire fleet using min/max expressions. While simple and accurate, this method is limited in capturing complex inter-dependencies and multi-dimensional constraints. Moreover, it doesn't scale as effectively for large systems or detailed optimization tasks.
%
%
In \cite{HARIGHI2024110732}, a stochastic optimization approach is used to develop an aggregation model that captures the flexibility potential of EV parking lots. The model estimates the parking lot's maximum flexibility margins (downstream and upstream) in advance, however, it does not capture the power boundaries and the minimum level capacities for the EVs.
A polytope-based aggregation approach was introduced in \cite{10121812}, offering accurate modeling of time-coupled and heterogeneous individual EV behaviors. However, the method is limited to aggregating EVs with identical arrival and departure times, making it unsuitable for capturing the diversity of charging demands or estimating flexibility across different times of the day. Additionally, it requires multiple iterations of the aggregation process, leading to higher computational costs for the aggregator.

Building on these estimation techniques, accurate forecasting of flexibility is crucial for real-time applications and operational planning.
%
In \cite{voss2018application}, forecasting models including ridge regression, K-nearest neighbors (kNN), Gaussian processes (GP), Support Vector Regression (SVR), and Artificial Neural Networks (ANN), were employed to improve the flexibility forecasting. However, due to the limitations of these models that can only learn a single target variable at a time, an iterative prediction approach was applied to generate forecasts over multiple time steps. 
Similarly, \cite{HUBER2020114525} employed techniques such as naive benchmarks, quantile regression, and kernel density estimation to forecast EV flexibility. In \cite{GENOV2024121969}, the accuracy and performance of Gaussian mixture models and LightGBM were assessed. However, these approaches remain limited as they typically rely on a single-output or cluster-based frameworks that do not address multivariate, multi-step interdependent forecasting as effectively as deep learning approaches.

In \cite{panda2024quantifyingaggregateflexibilityev}, SVR models were used for predicting aggregated EV flexibility. The two-step forecasting model predicts the aggregate charging profile within the flexibility request window, and the second step predicts the available flexibility in the same time window. 
Similarly, in \cite{HU20211101}, a combination of transformers and temporal convolution networks (TCNs) was used to forecast aggregated flexibility of an EV aggregator. Although the forecasting models in both papers were advanced and accurate, the underlying flexibility estimation strategy was restricted by its reliance on specific baselines and predefined flexibility services. Each unique combination of baseline and flexibility type (upstream or downstream) requires separate calculations and retraining, making it resource-intensive and less adaptable for varied scenarios. Moreover, the flexibility forecast becomes rigidly tied to predefined conditions, limiting its general applicability and scalability for dynamic energy market needs.

%

This paper presents novel methods for estimating and forecasting V2X flexibility potential of an EV fleet, based on an aggregate polytope representation. The methods are robust against individual uncertainties of EV owners behaviours as they are applied at an aggregate level, and  the reformulation of the V2X potential as a set of linear constraints allows the proposed methods to be integrated into different optimisation problems and therefore be applied for diverse V2X applications. We summarise the novelty of the paper as follow:
\begin{itemize}
    \item We propose a novel method for estimating V2X flexibility potential of an EV fleet based on a polytope representation,
    \item We leverage advanced multi-variate input multi-step multi-variate output  deep learning techniques for an accurate forecasting of V2X flexibility of an EV fleet, 
    \item We present examples of application of both V2X estimation and V2X forecasting to showcase the capabilities of both methods.
\end{itemize}
The paper is structured as follows: In Section II, a mathematical model is introduced to estimate V2X potential at an aggregate level. In Section III, state-of-the-art forecasting algorithms are leveraged to accurately predict V2X potential for future time horizons. 
The effectiveness of the method across diverse V2X applications is showcased in Section IV through a couple of case studies. Finally, Section V concludes the paper and discusses potential directions. 
\section{V2X flexibility estimation}
The operational constraints of a battery can be expressed in a polytope form \cite{polytope}: $\{p| Ap \leq b\}$, with $p$ is the power schedule of the battery and $A$ and $b$ are as follows:
\[A = 
\begin{bmatrix}
\textbf{I} & \textbf{0} \\
-\textbf{I} & \textbf{0}  \\
\textbf{0} & \textbf{I} \\
\textbf{0} & -\textbf{I} \\
(\eta_{in}\delta t) \Gamma & (\eta_{out}\delta t)\Gamma\\
-(\eta_{in}\delta t) \Gamma & -(\eta_{out}\delta t)\Gamma\\
\end{bmatrix}
\]
\[\Gamma =
\begin{bmatrix}
1 & 0 & 0 & \cdots & 0 \\
\alpha & 1 & 0 & \cdots & 0 \\
\alpha^2 & \alpha & 1 & \cdots & 0 \\
\vdots & \vdots & \vdots & \ddots & \vdots \\
\alpha^{T-1} & \alpha^{T-2} & \alpha^{T-3} & \cdots & 1 \\
\end{bmatrix}
\]
\[b=
\begin{bmatrix}
\textbf{P}_{max} \\
\textbf{0} \\
\textbf{0}  \\
-\textbf{P}_{min} \\
C_{max}-\alpha C_{0}  \\
\vdots \\
C_{max}-\alpha^TC_{0}  \\
\alpha C_{0} - C_{min}  \\
\vdots \\
\alpha^T C_{0} - C_{min}  \\
\end{bmatrix}
\]
The parameters figuring in the matrix $A$ and vector $b$ are the technical parameters of the battery and are summarized in Table \ref{tab:battery_params}.
\begin{table}
\caption{Description of the parameters in matrix  $A$ and vector $b$}
\label{tab:battery_params}
\begin{tabular}{|l| l|} 
\hline 
Parameter & Description \\
 \hline
 $\alpha$ & self-discharging rate \\
 $\eta_{in}$, $\eta_{out}$ & charging, discharging efficiency \\
 $\textbf{P}_{max}$, $\textbf{P}_{min}$ & maximum, minimum charging rate vectors \\
 $C_{max}$, $C_{min}$ & Maximum and minimum allowed energy levels \\
 $C_{0}$ & Initial energy level of the battery at time step $t=0$ \\
 \hline
 $C_{arr}$, $C_{dep}$ & Energy levels of the EV battery at arrival and  departure \\
 \hline
\end{tabular}
\end{table}
\textbf{I} and \textbf{0} are respectively the unity and zero vectors of size $T$, with $T$ being the time horizon for the battery operation, and $\delta t$ the time step of the operation. 
The first two blocks of the constraint $Ap\leq b$ define the upper and lower bounds on the charging rate of the battery. Similarly, the second two blocks define the upper and lower limits on the discharging rates, and the final two blocks limits the battery capacity in the range of its maximum and minimum capacity allowance. Hence, the size of A is (6$T$, 2$T$) and b is (6$T$,1). 

With relevant adjustments, the battery model can be adopted for EVs:
\begin{enumerate}
    \item The vector elements of $\textbf{P}_{max}$ and $\textbf{P}_{min}$ in the time slots where the EV is not available will be equal to zero.
    \item $C_0$ for a battery is equivalent to the capacity at arrival $C_{arr}$ for an EV.
    \item the last $T-t_{dep}$ rows of vector $b$ can be modified to $\alpha^{t_{dep}+i} C_{arr} - C_{dep}$ to force the battery to be charged at least at the desired level $C_{dep}$ at time of departure $t_{dep}$.
    \item The first [0-t$_{arr}$] elements of the minimum capacity limit block are set to 0. 
\end{enumerate} 
With these cited adjustments, vector $b$ becomes:
\[b_{ev}=
\begin{bmatrix}
\textbf{P}_{max} \\
\textbf{0} \\
\textbf{0}  \\
-\textbf{P}_{min} \\
\hdashline 
0\\
\vdots \\
0 \\
C_{max} - \alpha C_{arr}  \\
\vdots \\
C_{max} - \alpha^{t_{dep} -t_{arr}} C_{arr} \\
\vdots \\
C_{max} - \alpha^{t_{dep} -t_{arr}} C_{arr} \\
\hdashline 
0\\
\vdots \\
0 \\
\alpha C_{arr} - C_{min}  \\
\vdots \\
\alpha^{t_{dep} -t_{arr}} C_{arr} - C_{dep}\\
\vdots \\
\alpha^{t_{dep} -t_{arr}} C_{arr} - C_{dep}  \\
\end{bmatrix}
\]
To enhance readability, horizontal lines were added to separate the power bounds block from the maximum capacity limit block and the minimum capacity limit block.
Following the demonstration in \cite{polytope}, an aggregation of EVs $(EV_i)_{i \in [1, N]}$ can also be represented as a polytope: $$\{p_{agg}| A_{agg}p_{agg} \leq b_{agg}\}$$
Assuming all EVs in the aggregation have the same self discharge rate $\alpha$ and the same charging/discharging efficiencies $\eta_{in}$/$\eta_{out}$, matrix $A$ will be the same for all EVs, and the aggregate polytope can be expressed as: 
\begin{equation}
    P = \{p_{agg}| A_{agg}p_{agg} \leq b_{agg}\}
    \label{eq:p_agg}
\end{equation}
with:
\[
\left\{
\begin{aligned}
    A_{agg} &= A = A_1= A_2 =  \cdots = A_N \\
    p_{agg} &= \sum_{i=1}^{N} \textbf{p}_i\\
    b_{agg} &= \sum_{i=1}^{N} b_i
\end{aligned}
\right.
\]
In a compact form $b_{agg}$ can be written as:
\[b_{agg}=
\begin{bmatrix}
\textbf{P}^{agg}_{max} \\
\textbf{0} \\
\textbf{0}  \\
-\textbf{P}^{agg}_{min} \\
\textbf{C}^{agg}_{max} \\
\textbf{C}^{agg}_{min} \\
\end{bmatrix}
=
\begin{bmatrix}
\sum_{i=1}^{N} \textbf{P}^i_{max} \\
\textbf{0} \\
\textbf{0}  \\
-\sum_{i=1}^{N} \textbf{P}^i_{min} \\
\sum_{i=1}^{N} \textbf{C}'^{,i}_{max} \\
\sum_{i=1}^{N} \textbf{C}'^{,i}_{min} \\
\end{bmatrix}
\]
with:
\[\textbf{C}'^{,i}_{max}=
\begin{bmatrix}
0\\
\vdots \\
0 \\
C^i_{max} - \alpha C^i_{arr}  \\
\vdots \\
C^i_{max} - \alpha^{t^i_{dep} -t^i_{arr}} C^i_{arr} \\
\vdots \\
C^i_{max} - \alpha^{t^i_{dep} -t^i_{arr}} C^i_{arr} \\
\end{bmatrix}
\]
and
\[\textbf{C}'^{,i}_{min}=
\begin{bmatrix}
0 \\
\vdots \\
0
\\
\alpha C^i_{arr} - C^i_{min} \\
\vdots \\
\alpha^{t^i_{dep}-t^i_{arr}} C^i_{arr} - C^i_{dep} \\
\vdots \\
\alpha^{t^i_{dep}-t^i_{arr}} C^i_{arr} - C^i_{dep} \\
\end{bmatrix}
\]

We note that the assumption of similar self discharge rate and charging/discharging efficiencies is a reasonable approximation for the following reasons:
\begin{itemize}
    \item Aggregation smooths individual variabilities: variations in individual EV parameters tend to cancel each other out, making the collective behavior more stable and less sensitive to minor differences between individual vehicles,
    \item Standardization of EV battery models: Advances in battery technology have led to convergence in the efficiency of charging and discharging processes across different EVs, reducing the magnitude of variability. Many EVs are composed of standardized battery models with comparable parameters, making the considered assumptions realistic and consistent with industry trends. 
\end{itemize}
At the bulk system or transmission system level, this method offers several significant benefits. By capturing the collective behavior of a large fleet of EVs, the aggregation approach enables the system operator to treat the entire fleet as a single, flexible resource that can be utilized to support grid stability and efficiency. The aggregated model allows for the incorporation of the fleet's flexibility into system-wide optimization problems, such as those related to balancing supply and demand, managing congestion, or integrating renewable energy sources. This enhances the ability to respond to market signals or grid emergencies with greater precision, offering a scalable solution that can be integrated into broader grid management strategies. Similarly, if we consider a charging station with N EVs, the set of linear constraints $A_{agg}p_{agg} \leq b_{agg}$ captures all the technical requirements of the EV fleet over the considered time horizon $T$, and can therefore be included in any optimisation problem to schedule the charging/discharging operations of the charging station in response to a market signal or an incentive while respecting the operational constraints of the EV batteries. In other words, solving an optimisation problem with system (\ref{eq:p_agg}) provides the charging station operator with insights into the potential discharging capacity that can be deployed to the grid at each time step $t$, while ensuring compliance with EVs constraints, and potentially other constraints embedded within the optimization problem (e.g., network constraints).
We note that solving a problem that incorporates the aggregate polytope constraints will return an aggregate solution $p_{agg}$ at an aggregate level. The formula that recovers a feasible solution $\textbf{p}$ from an approximate aggregate solution $p_{agg}$ is presented in \cite{polytope} and is as follows:
\begin{equation}
\begin{aligned}
\min_{p_{i}, i=1, \dots, N} \quad & \left\Vert p_{agg} - \sum_{i=1}^{N}p_i \right\Vert  \\
\textrm{s.t.} \quad & A_ip_i \leq b_i, \quad i=1, \dots, N
\end{aligned}
\end{equation}
\subsection{Example}
\begin{figure}
    \centering
    \includegraphics[width=0.5\textwidth]{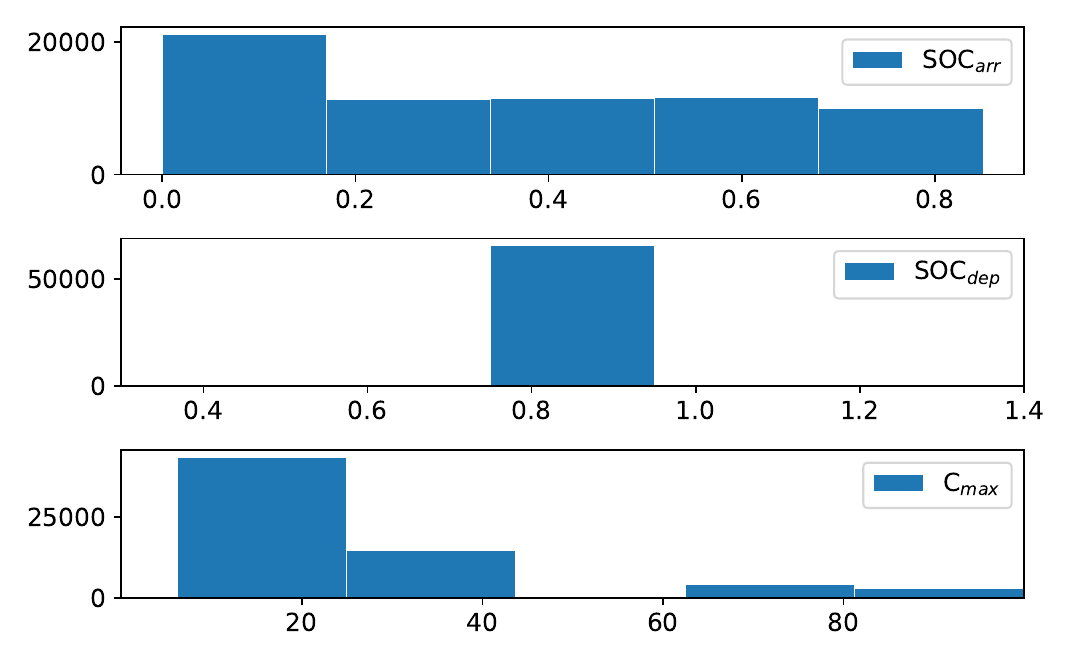}
    \caption{\footnotesize Overview of the EVs statistics for the considered dataset: The majority of the cars come with a SOC$_{arr}$  less than 0.2 and leave with a SOC greater than 0.8. The battery capacities are concentrated around 20 kWh with few of them in the range of 40, 80 and 100 kWh.}
    \label{fig:stats}
\end{figure}

\begin{figure}
    \centering
    \includegraphics[width=0.5\textwidth]{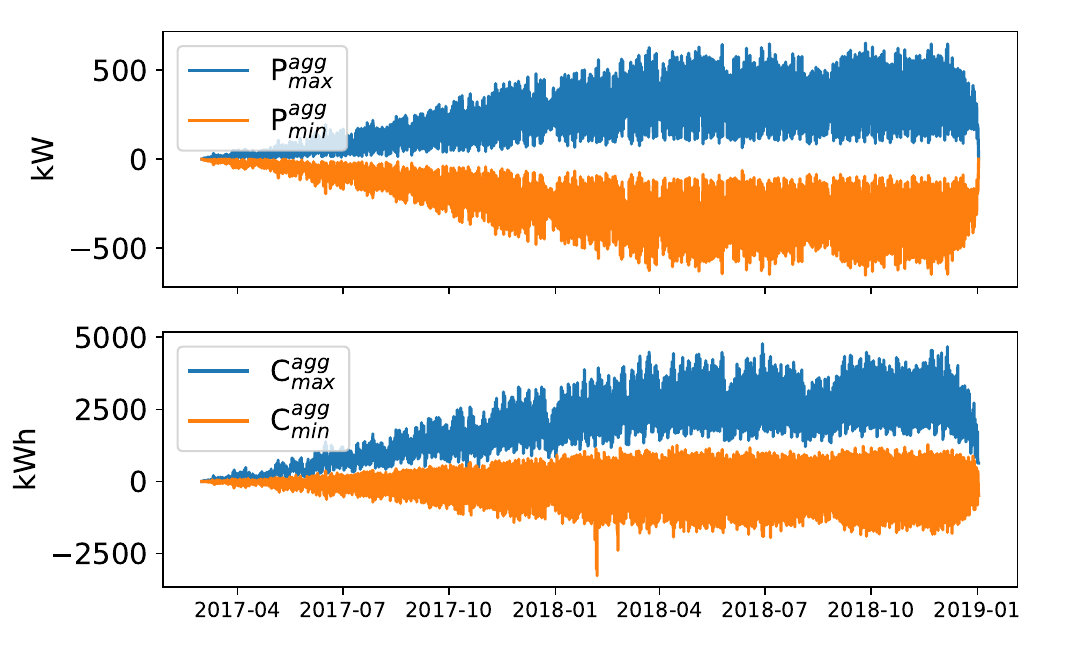}
    \caption{\footnotesize The vectors of the constructed $b_{agg}$ from the considered dataset. $P_{min}^{agg}$ is equivalent to the negative value of $P_{max}^{agg}$ as we considered that the EVs have symmetrical charging/discharging capacities. Roughly, the aggregator can be seen as an EV with 500 kW maximum charging rate and 4000 maximum capacity.}
    \label{fig:bagg}
\end{figure}

We work with the \textit{Crowd Charge} dataset \cite{EV-dataset}, 
a dataset that collected charging session data from approximately 700 EV owners who participated in an 18-month trial for a home smart charging initiative conducted in the UK. 
This dataset was selected because it is the only open-source dataset that has recorded the capacity of the tracked EV batteries $C_{max}$, essential for the construction of the individual $b_{ev}$ vectors. The capacity at arrival $C_{arr}$ and the desired capacity at departure $C_{dep}$ are not available in the dataset and were approximated as shown in Eqs. (\ref{eq:B8}), with $C_{cons}$ (kWh) is the consumed energy during the session. 
We select a horizon of one day with 15 min resolution.
Since $b_{ev}$ is formulated over a finite horizon $T$ (i.e., one day in our case) and the EV availability can span across multiple days, particularly in the case of home charging, we supposed that the battery capacity at the end of a day should be proportional to the final desired capacity. For example, if an EV arrives with a capacity $C_{arr}$, wants to leave with a capacity $C_{des}$, and is parking for 6 time slots, with 4 time slots on day 1 and 2 time slots in day 2, then  $C_{arr}^{day 1} = C_{arr}$ , $C_{des}^{day 1} = \frac{4}{6}C_{des}$, $C_{arr}^{day 2} = C_{des}^{day 1}$ and $C_{des}^{day 2} = C_{des}$.
\begin{table*}
\centering
\begin{minipage}{1\textwidth}
 \begin{equation}
 \left\{ 
\begin{aligned}
\label{eq:B8}
C_{arr} = 0.2C_{max} \quad & \mbox{and} & C_{dep} = C_{arr} + C_{cons} \quad & \mbox{if} & 0.2C_{max}  + C_{cons} \leq C_{max}     \\
C_{arr} = C_{max} - C_{cons} \quad & \mbox{and} & C_{dep} = C_{max}  \quad &  \mbox{if} & 0.2C_{max}  + C_{cons} > C_{max}    \\
C_{min} = 0 & &
\end{aligned}
\right.
\end{equation}

\medskip
\hrule
\end{minipage}
\end{table*}
Figure \ref{fig:stats} gives an overview of the different EV batteries information following the application of the mentioned approximations. Most EVs arrive with a state of charge $SOC_{arr}$ around 0-0.2\%, and depart with a $SOC_{dep}$ of 0.8\% or larger. The EVs that participated in the trial have capacities ranging between 20 kWh to 100 kWh, with the majority concentrated around the 20 kWh.
The $b_{ev}$ vectors are then constructed using OPLEM platform \cite{oplem} and summed up according to Eq. \ref{eq:p_agg} to return $b_{agg}$. The source code for computing $b_{agg}$ can be found in the Github repository of the project\footnote{\url{https://github.com/GridLab-NTU/V2X-prediction}}. Figure \ref{fig:bagg} displays the aggregate maximum charging/discharging rates and the aggregate maximum/minimum capacities for the duration of the project trial between 01 March 2017 and 01 January 2019. Roughly speaking, an aggregator managing this fleet of EVs can be seen as an EV with 500 kW maximum charging rate, 4000 maximum capacity and around 2000 minimum capacity.
An example demonstrating how can we leverage such information to gain insights is presented in the first case study of Section \ref{sec:cs}.
\section{V2X flexibility forecasting}
In real-world applications, transmission system operators and EV aggregators such as charging station operators (CSOs) must assess the availability and demand 
of EVs ahead of time to effectively allocate resources. As the vector $b_{agg}$ captures all the technical constraints of an aggregate EV fleet, we assess in this paper the capability of forecasting $b_{agg}$ using advanced forecasting models.
\subsection{Setting}
We use the same dataset as in the previous section and train different deep learning forecasting models to predict the day-ahead values of $\textbf{P}_{max}^{agg}$, $\textbf{C}_{max}^{agg}$ and $\textbf{C}_{min}^{agg}$  vectors. $\textbf{P}_{min}^{agg}$ was not forecasted as it is equal to the negative value of $\textbf{P}_{max}^{agg}$. For EV fleet containing EVs with asymmetric charging/discharging rates, $\textbf{P}_{min}^{agg}$ should be also predicted. 
The training was conducted in two phases. In the first phase, we explored various deep learning (DL) architectures, including RNN, LSTM, bi-LSTM, and GRU, to identify the configurations that consistently delivered the best performance across different settings. The details of the explored architectures can be found in the Appendix \ref{app1}, Table \ref{tab:train_1}.
Based on the results from this initial round, we selected attention-based CNN-LSTM and Transformer models for further evaluation. In the second phase, we performed hyperparameter tuning on these selected models to enhance their accuracy and optimize their predictive capabilities. The specifics of the models and the hyperparameters adjusted during tuning are detailed in the Appendix \ref{app1}, Table \ref{tab:train_2}.

\begin{figure}
    \centering
    \includegraphics[width=0.45\textwidth]{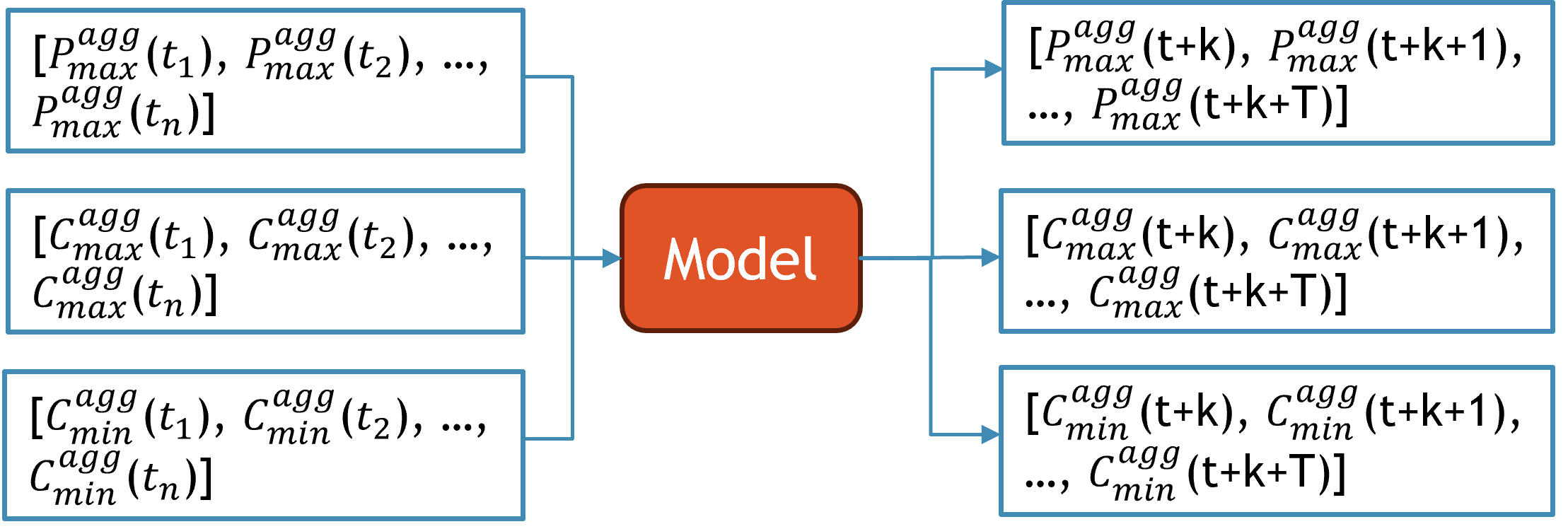}
    \caption{\footnotesize Multivariate input multi-step multivariate output with n historical data, $k$ lead time and $T$ steps.}
    \label{fig:model-structure}
\end{figure}

Our objective is to design a model capable of simultaneously forecasting three variables over a defined horizon $T$. Fig. \ref{fig:model-structure} illustrates the general structure of this model. This model input consists of several vector variables (e.g., $\textbf{P}_{max}^{agg}$, $\textbf{C}_{max}^{agg}$ and $\textbf{C}_{min}^{agg}$ in our case), each with $n$ historical values. The model then forecasts multiple vectors as the output, where each output vector predicts a sequence of $T$ future time steps over a lead time of $k$ steps into the future.
Unlike univariate output models, which require a separate model for each variable, this approach consolidates the forecasting process into a single model for all three variables, leveraging the relationships across both the different variables and the multiple time steps. However, this integrated approach presents a greater challenge in training to achieve satisfactory accuracy.

The input data included values of the three main features: $\textbf{P}_{max}^{agg}$, $\textbf{C}_{max}^{agg}$ and $\textbf{C}_{min}^{agg}$, 
and is structured as follows:
\begin{itemize}
    \item Current Value: The current values of the 3 variables $\textbf{P}_{max}^{agg}$, $\textbf{C}_{max}^{agg}$ and $\textbf{C}_{min}^{agg}$, namely, $P_{max}^{agg}(w,d,t)$, $C_{max}^{agg}(w,d,t)$ and $C_{min}^{agg}(w,d,t)$, with $w$ referring to the index of the current week, $d$ the current day and $t$ to the current time slot.
    \item Quarter-hourly History: Values of the variables for the past 23 hours, capturing short-term dependencies. For $\textbf{P}_{max}^{agg}$ for instance, we have $P_{max}^{agg}(w,d, t-\delta t)$, ..., $P_{max}^{agg}(w,d, t-T)$, 
    \item Daily History: Values of the variables for the same time across the past 6 days, capturing daily seasonality, e.g., $P_{max}^{agg}(w,d-1, t+k)$, ..., $P_{max}^{agg}(w,d-6, t+k)$, 
    \item Weekly History: Values of the variables for the same time across the past 5 weeks, capturing weekly patterns, e.g., $P_{max}^{agg}(w-1,d, t+k)$, ..., $P_{max}^{agg}(w-5,d, t+k)$.
    \item Target: The models predict the 96 values of the three variables beginning from the $k^{th}$ time slot ahead of the current moment, e.g., $P_{max}^{agg}(w,d, t+k)$, ..., $P_{max}^{agg}(w,d, t+k+96)$
\end{itemize}
We utilized Python, TensorFlow (DL models), and Scikit-learn (for baseline model) for model development and training. The models were trained using Google Colab, a cloud-based platform with NVIDIA Tesla T4 GPUs and 12 GB of RAM resources to accelerate the training process. However, the 12 GB of RAM was not consistently utilized during all model training processes or for every model iteration. A (90\%-10\%) rate was selected for the train-test split. During the training, we tried several regularization techniques, such as Batch Normalization, L1, and L2 regularization, however the models were overfitted, likely due to the simplicity of our dataset. Given this, we trained the models without regularisation, but we recommend including a regularisation technique for richer datasets.
\subsection{Results}
\begin{figure}
    \centering
        \begin{subfigure}{\linewidth}
            \centering
            \includegraphics[width=\textwidth]{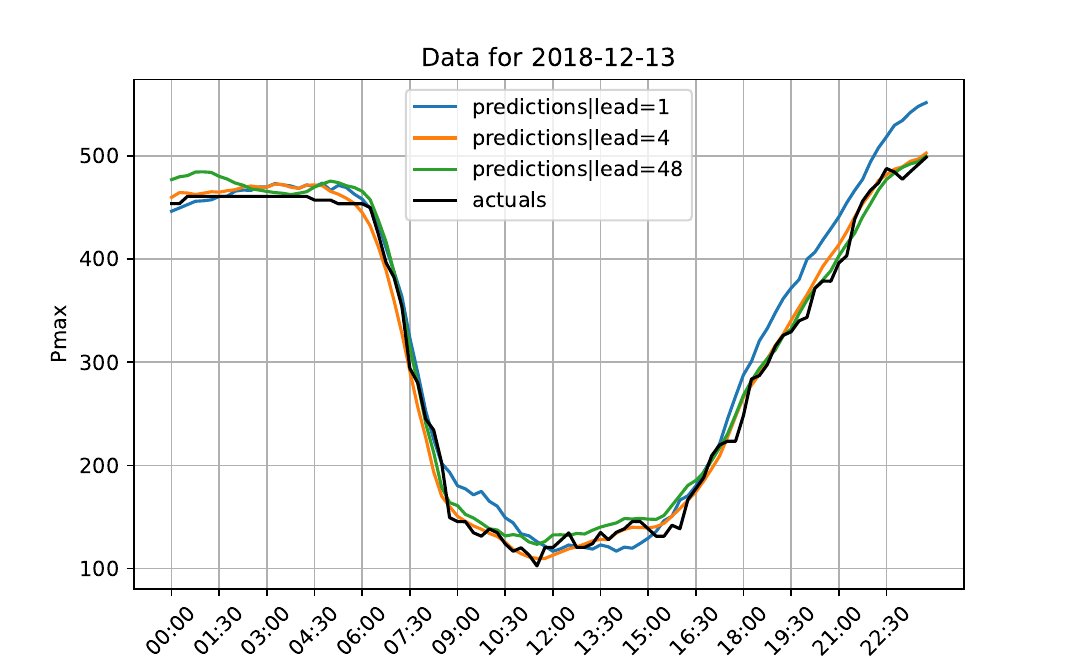}
            \caption{$P_{max}$ forecast}
        \end{subfigure}
        \hfill
        \begin{subfigure}{\linewidth}
            \centering
            \includegraphics[width=\textwidth]{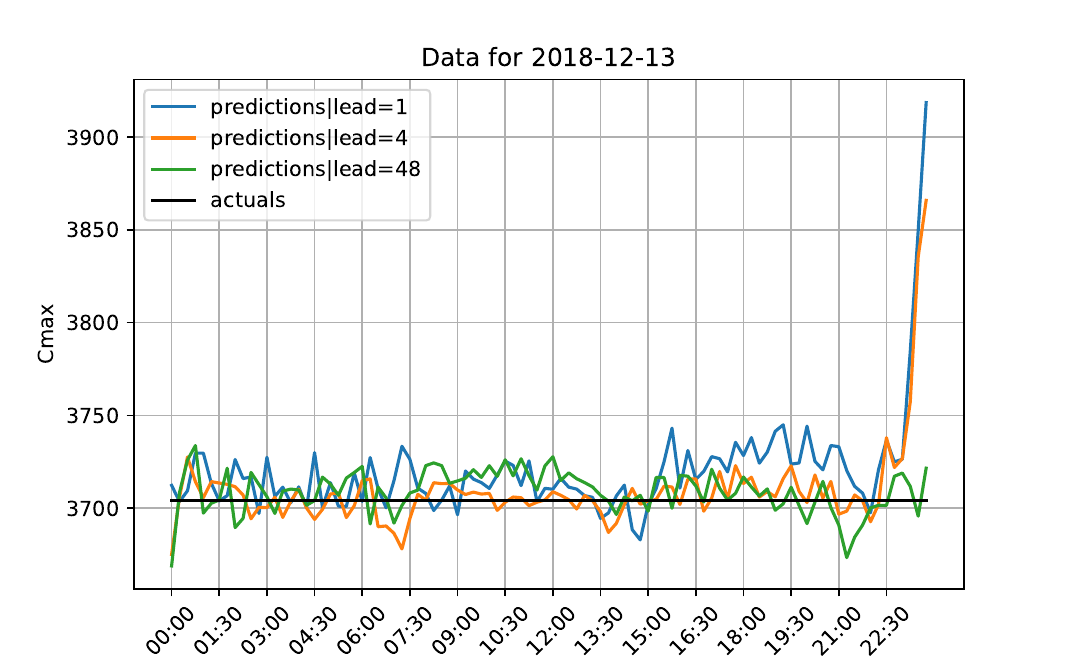}
            \caption{$C_{max}$ forecast}
        \end{subfigure}
        \hfill
        \begin{subfigure}{\linewidth}
            \centering
            \includegraphics[width=\textwidth]{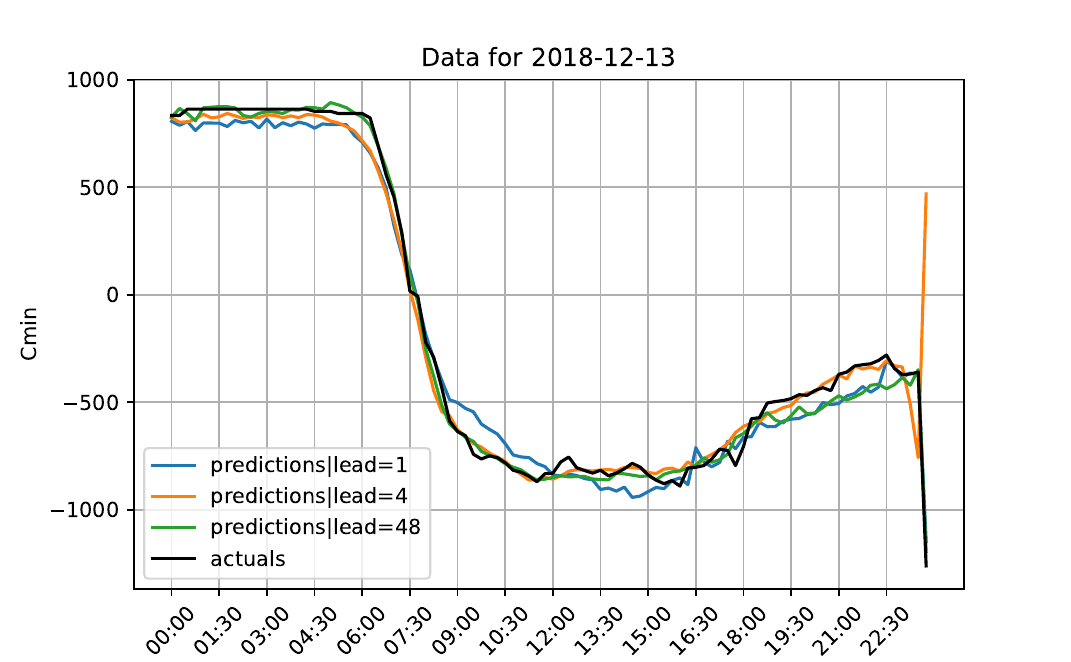}
            \caption{$C_{min}$ forecast}
        \end{subfigure}
    \caption{The forecast of $b_{agg}$ vectors for one day.}
    \label{fig:bagg_fore}
\end{figure}

Three look-ahead values ($k$) were determined for the forecasting:
\begin{itemize}
    \item $k=$1: model trains to forecast the value of the $b_{agg}$ vector for 96 time slots (equivalent of 1 day) beginning from the next time slot. As the time resolution is 15 min, one step-ahead forecast returns one day forecast of the vector $b_{agg}$ starting from the next 15 min slot. For example, at time step 36 (9~am), the model forecasts $b_{agg}$ from 9am15 until 9~am00 of the next day. This type of forecasts is beneficial for aggregators who want to participate in real or near-real time markets such as FCR, aFRR and mFRR.
    \item $k=$4: model trains to forecast $b_{agg}$ for one day beginning from the 4$^{th}$ next time slot (i.e., hour-ahead). For instance, at 9~am the model forecasts $b_{agg}$ from 10~am to 9~am45 of the next day. This foresight not only helps aggregators capitalize on intra-day trading opportunities, optimize resource utilization, and avoid potential penalties for imbalances—thereby improving profitability and operational efficiency—but also provides critical insights for bulk system and transmission system operators. Accurate 1-hour ahead forecasts enable these operators to maintain grid stability, efficiently manage load balancing, and anticipate fluctuations in supply and demand.
    \item $k=$48: model trains to forecast $b_{agg}$ for one day beginning from the 48$^{th}$ next time slot (i.e., half-day-ahead). At 12~pm, the model forecasts $b_{agg}$ for the next day from 00~am00 to 11~pm45. When engaging in day-ahead markets, half-day ahead forecasts are essential for strategic planning and resource allocation. These forecasts provide the necessary information to make informed decisions about energy bids and schedules for the next day.
\end{itemize}
The parameters of the top two performing Transformer and attention-based CNN-LSTM models are presented in Table \ref{tab:winning-architectures}, and the forecasting errors (RMSE) are presented in Table \ref{tab:rmse}. The RMSE results are shown alongside a RandomForest (RF) model, which serves as the baseline. Note that the baseline model is not capable of producing multi-variate, multi-step outputs; therefore, three separate models were trained, each corresponding to a different variable.
\begin{table}[h]
    \centering
    \begin{tabularx}{.45\textwidth}{|X|X|l|l|l|}
    \hline
    Model & Parameter & lead=1 & lead=4 & lead=48 \\
    \hline
    Transformer & N. of layers & 4 & 5 & 5  \\
    \cline{2-5}
                & d\_model & 128 & 128 & 128 \\
    \cline{2-5}
                & N. of heads &  8 & 16 & 16 \\
    \cline{2-5}
                & dff & 512 & 128 & 128  \\
    \hline
    Attention-based CNN-LSTM & N. of CNN layers & 1 & 2 & 1  \\
    \cline{2-5}
                             & Filters & 64 & 64-128 & 128 \\
    \cline{2-5} 
                & N. of LSTM layers & 2 & 3 & 3 \\
    \cline{2-5}             
                 &  units & 64-128 & 128-128-64 & 128-128 \\
    \hline
    \end{tabularx}
    \caption{Parameters of the most performing DL architectures}
    \label{tab:winning-architectures}
\end{table}
\begin{table}[h]
    \centering
    \caption{Forecasting results.}
    \label{tab:rmse}
    \begin{tabular}{|c|c|c|c|c|}
    \cline{3-5}
      \multicolumn{1}{c}{}&\multicolumn{1}{c|}{} & RF & CNN-LSTM & Transformer\\
    \hline
       \multirow{4}{*}{$k=1$}  & $P_{max}^{agg}$ & 32.31 & 12.83    &  10.96       \\
                              & $C_{max}^{agg}$ & 39.84 & 58.20    &  35.37       \\
                              & $C_{min}^{agg}$ & 209.10 & 90.01    &  60.35       \\
    \cline{2-5}
                              & Average         & 93.75 & 53.68    & 35.56          \\
    \hline
       \multirow{4}{*}{$k=4$}  & $P_{max}^{agg}$ & 32.7 & 11.62    &  9.00      \\
                              & $C_{max}^{agg}$ & 36.06 & 52.65    &  30.72      \\
                              & $C_{min}^{agg}$ & 211.04 & 106.02   &  48.77       \\
    \cline{2-5}
                              & Average         & 93.26 & 56.61    & 29.5          \\
    \hline
     \multirow{4}{*}{$k=48$}  & $P_{max}^{agg}$ & 20.6 & 11.35    &  9.30      \\
                              & $C_{max}^{agg}$ & 36.06 & 52.36    &  33.72      \\
                              & $C_{min}^{agg}$ & 211.04 & 69.06   &  50.01       \\
    \cline{2-5}
                              & Average         & 88.96 &  44.25   &  31.01       \\
    \hline
\end{tabular}
\end{table}
The results are promising as they demonstrate that using these forecasting algorithms, we were able to reach an acceptable accuracy, less than 6.59\% error forecast.
Figure \ref{fig:bagg_fore} shows the results of the prediction for one day selected from the test period data.

\section{Case studies}
\label{sec:cs}
We demonstrate 2 scenarios on leveraging $b_{agg}$ to gain insights and participate in local energy markets and local flexibility markets.
We work with the same dataset as previously and suppose that the aggregator possesses historical data of charging sessions that have relevant information:
\begin{itemize}
    \item Charging/discharging rates of EV batteries,
    \item Connection/disconnection times of EV batteries,
    \item State of charge at arrival and maximum capacity of EV batteries,
    \item And optionally the desired state of charge at departure.
\end{itemize}
\subsection{Case study 1: Participation in long-term flexibility contracts}
\textit{Sustain} flexibility service, as defined by Open Networks \cite{ena}, entails the distribution system operator (DSO) to enter in a long-term contracts with the flexibility service provider (FSP) in which the FSP agrees to alter its power during a pre-defined window. This type of services presents a cost-effective alternative to network reinforcement and is being launched by some DSOs \cite{ssen, ukpn} in the UK. Call of tenders are submitted and FSPs send their bids composed of the pair \{flexibility, price\}.

In our case study, we investigate the potential of the aggregator to engage in a long-term upward flexibility contract with the DSO, and present a method that computes the amount of flexibility the aggregaror can bid as a response to a DSO call using the $b_{agg}$ vector.

First, the flexibility is measured as a deviation from a baseline schedule. Here, we assume that the aggregator baseline is an optimal schedule of its resources in response to a day-ahead market price signal: 
\begin{equation}
\begin{aligned}
\underset{(p_{agg})}{min} \quad & \sum\limits_{t \in \mathcal{T}}   \lambda^t_{imp} p^t_{agg,ch} - \lambda^t_{exp} p^t_{agg,dis}
& \\
\textrm{s.t.} \quad &  A_{agg}p_{agg} -b_{agg} \leq 0 
\end{aligned}
\label{eq: baseline}
\end{equation} 
where $\lambda^t_{imp}/\lambda^t_{exp}$ are the day-ahead import/export prices and $p_{agg}=[p^t_{agg,ch}, p^t_{agg,dis}]^T$ is the charging/discharging aggregated power.

\begin{figure*}
    \centering
        \begin{subfigure}{.49\linewidth}
            \centering
            \includegraphics[width=\textwidth]{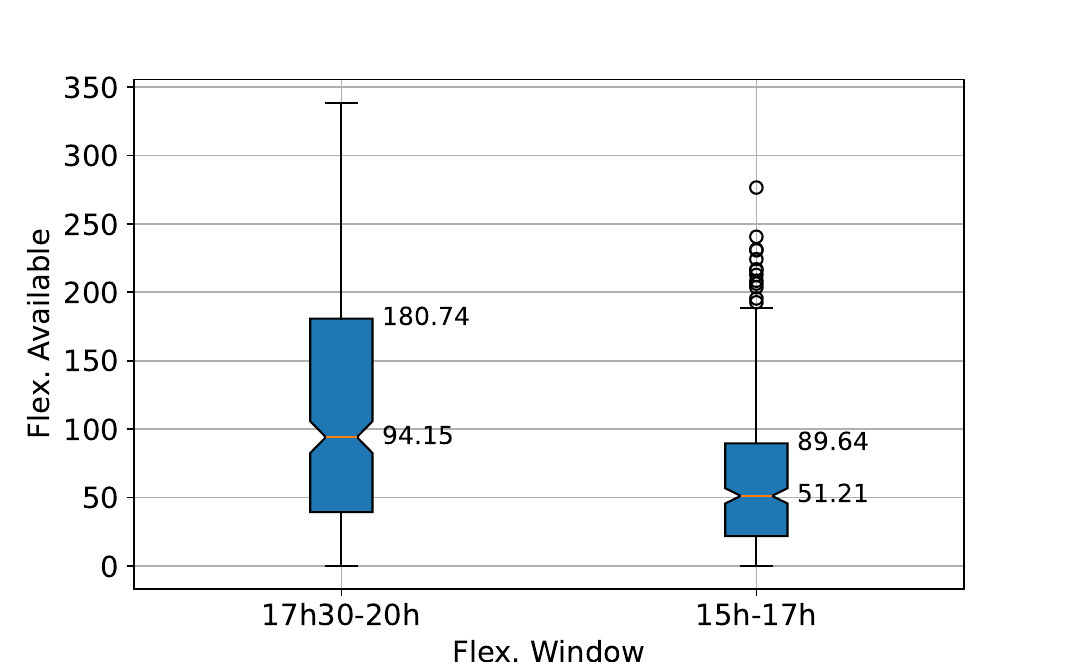}
            \caption{The period of 01 March 2017 to 28 February 2018.}
            \label{fig:flex-stats}
        \end{subfigure}
        \begin{subfigure}{.49\linewidth}
            \centering \includegraphics[width=\textwidth]{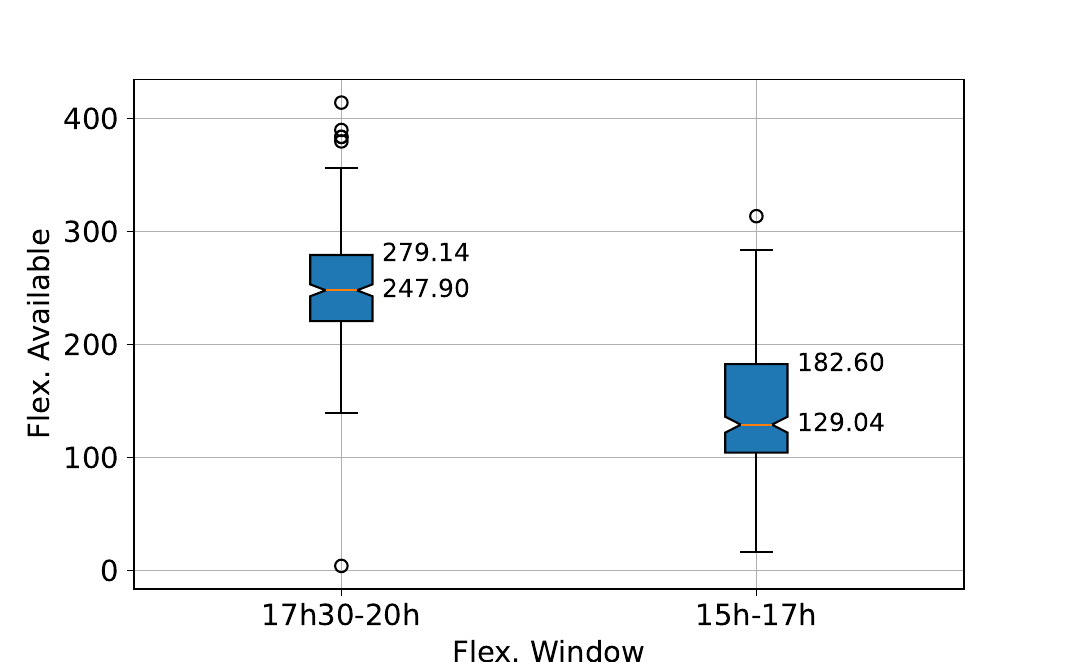}
            \caption{The period of 01 March 2018 to 31 December 2018.}
            \label{fig:flex-stats-period2}
        \end{subfigure}
\caption{Maximum flexibility provision for for two flexibility windows, in the period of feasibility study (left) and period of flexibility activation (right).}
\end{figure*}

Let us note $p_{agg}^{sch}$ the optimal schedule of the aggregator in response to the market signal, e.g., the result of the optimisation problem (\ref{eq: baseline}) and $[T_s, T_e]$ the time window of providing flexibility. The computation of the maximum upward flexibility the aggregator can offer in the time window $[T_s, T_e]$ can be retrieved using the following optimisation problem: 
\begin{equation}
\begin{aligned}
\underset{(flex, p_{agg})}{max} \quad & flex
& \\
\textrm{s.t.} \quad &  A_{agg}p_{agg} -b_{agg} \leq 0  \\
                    & p_{agg}^{sch} - p_{agg} \geq flex, \quad \forall t \in [T_s, T_e]
\end{aligned}
\label{eq: flex}
\end{equation} 
Running this procedure (i.e., problem (\ref{eq: baseline}) followed by problem (\ref{eq: flex})) for available historical data will provide the aggregator with insights into its flexibility capability.
Fig. \ref{fig:flex-stats} shows the results of running the procedure for days between 01 March 2017 and 28 February 2018. We find that the aggregator can provide more upward flexibility in the time window 17h30-20h than in 15h-17h, which is expected as this aggregator is managing home-level charging with more availability in the late afternoon than early afternoon.
We suppose two scenarios, in one of them the aggregator bids the median value of the flexibility and in the other it bids the third-quantile value. We then assess the reliability of the aggregator to deliver the promised flexibility over the period of 01 March 2018 to 31 December 2018. Table \ref{tab:success-rates} shows the results for both scenarios. Successful delivery refers to the cases where the aggregator was able to deliver 90\%$^+$ of its bid, partial delivery to [50, 90[\% of the bid, and failure of delivery when the delivered flexibility is less than 50\% of the agreed capacity.
\begin{table*}
    \centering
    \caption{Success rates for flexibility delivery over the period of 01 March 2018 to 31 December 2018.}
    \label{tab:success-rates}
    \begin{tabular}{c|ccc|ccc|}
    \cline{2-7}
         & \multicolumn{3}{c}{17h30-20h} & \multicolumn{3}{|c|}{15h-17h}  \\
          & Success & Partial success & Failure & Success & Partial success & Failure \\
    \hline
    Median & 99.68\% & 0.0\%   & 0.32\% & 99.67\% & 0.0\% & 0.33\%\\
    Q3     & 99.0\% & 0.65\% & 0.32\% & 92.81\% & 6.86\% & 0.32\% \\
    \hline
    \end{tabular}
\end{table*}
Whether bidding the median or the 3$^{rd}$ quantile value, the aggregator was 99\% of the time successful in delivering the bid flexibility. 
Upon investigating the underlying reasons and closely examining the available flexibility during this second period, as shown in Fig. \ref{fig:flex-stats-period2}, we observed a notable increase in the aggregator's flexibility potential. This can be attributed to the inclusion of new EV members joining the project trials, which explains the high success rate of flexibility delivery.
\subsection{Case study 2: Day ahead scheduling of the charging station operation in response to dynamic operating envelop signal}
Dynamic operating envelops (DOEs) define the real-time limits for energy export and import at a connection point \cite{pesgm-doe}. Unlike static limits, DOEs adjust dynamically based on factors such as grid capacity, network conditions, and local demand. This approach enables the seamless integration of renewable energy sources and optimizes grid efficiency.

In this case study, we consider a scenario where the DSO sends day-ahead DOE signals to flexible customers, such as EV aggregators, who adjust their next-day resource scheduling accordingly. 
The mathematical formulation of this problem involves minimizing the energy bill, subject to the operational constraints of the aggregator and the DOE limits:
\ \begin{equation}
\begin{aligned}
\underset{p_{agg}}{min} \quad & \sum\limits_{t \in \mathcal{T}}   \lambda^t_{imp} p^t_{agg,ch} + \lambda^t_{exp} p^t_{agg,dis}
& \\
\textrm{s.t.} \quad &  A_{agg}p_{agg} -b_{agg} \leq 0 \\
                    & p^t_{agg,ch} \leq p^t_{doe, imp}, \quad \forall t\in \mathcal{T}\\
                    & p^t_{agg,dis} \geq p_{doe, exp}, \quad \forall t\in \mathcal{T}\\
\end{aligned}
\label{eq: baseline2}
\end{equation} 
Where $p_{agg} = [p^t_{agg,ch}, p^t_{agg,dis}]^T$ is the scheduled power output, $p^t_{doe, imp}, p^t_{doe, exp}$ are the import/export DOE limits and $\lambda^t_{imp}, \lambda^t_{exp}$ are the day-ahead  import/export prices.

\begin{figure}
    \centering
    \includegraphics[width=0.45\textwidth]{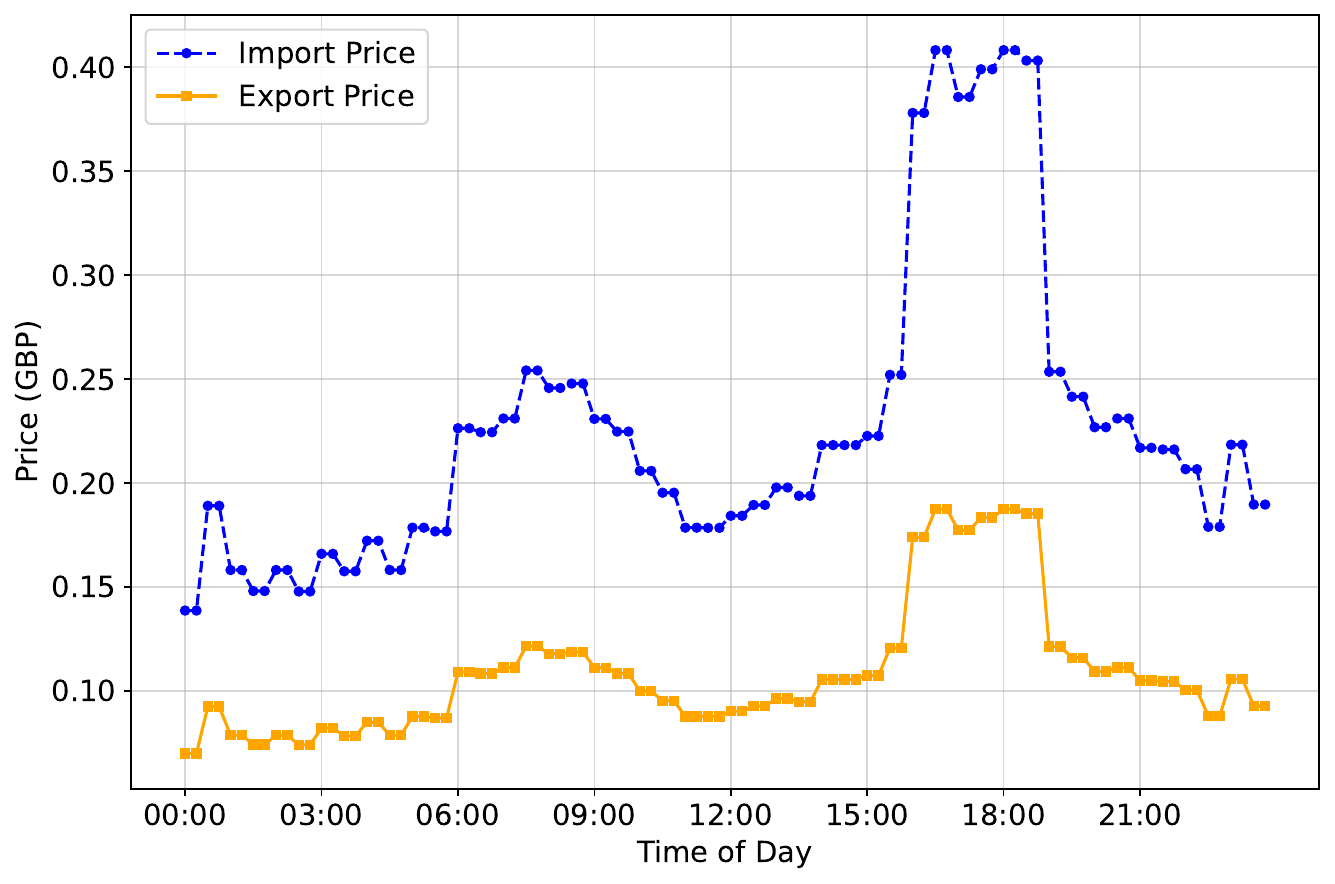}
    \caption{Import and export energy prices.}
    \label{fig:Tariffs}
\end{figure}

\begin{figure}
    \centering
    \includegraphics[width=\linewidth]{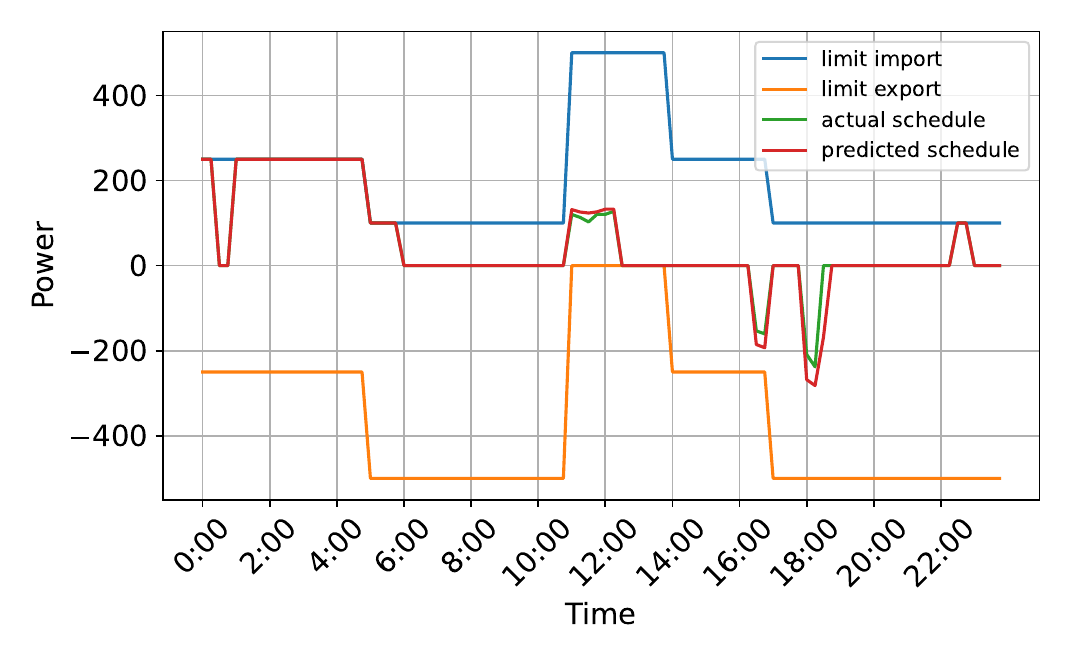}
    \caption{Aggregator schedule (actual and predicted) in response to DOE signal.}
    \label{fig:CS2}
\end{figure}

We performed the optimization (\ref{eq: baseline2}) for both the forecasted $b^{fore}_{agg}$ vector and the actual $b^{act}_{agg}$ vector, resulting in two types of power schedules one using the predicted $b^{fore}_{agg}$ values and one using the $b^{act}_{agg}$ values. 
The time-varying tariffs \textit{Agile Octopus April 2024 V1} (for imports) and \textit{Agile Outgoing Octopus May 2019} (for exports) were utilized, with the tariffs illustrated in Fig. \ref{fig:Tariffs} for November 20, 2024.
The selected DOEs for the simulation reflect a system with high solar energy penetration, imposing high (low) import (export) limits around midday, moderate import/export limits in the morning and around late afternoon, and a low (high) import (export) limits in the evenings. 
The results in Fig. \ref{fig:CS2} demonstrate the effectiveness of our prediction method for the $b$ vector, showing a slight difference between the two schedules. 
%
\section{Summary and conclusions}
In this paper, we presented a comprehensive approach for estimating and forecasting V2X flexibility. Our proposed method accurately estimate V2X flexibility potential at an aggregate level and utilize DL architectures to predict future flexibility with high precision. Through two case studies, we demonstrated the practical applicability and effectiveness of our methods in real-world scenarios.

However, while our approach shows promise, several limitations must be addressed to improve the reliability and robustness of V2X flexibility forecasts. A significant challenge is the lack of detailed data, particularly regarding electric vehicle (EV) battery capacity and state of charge at arrival (SOC$_{arr}$). Many charging stations currently do not retrieve these critical pieces of information during charging sessions, which are essential for generating reliable V2X flexibility forecasts. To overcome this limitation, there is a pressing need for regulations that mandate the collection of such data and for the development of technical equipment capable of retrieving it. 

For future work, it is crucial to collect larger datasets with comprehensive and relevant information. This will not only enable more efficient training of DL methods but also enhance the robustness and accuracy of the forecasts, ultimately contributing to more reliable and actionable insights for energy network operators and other stakeholders in the V2X ecosystem.



\appendices
\section{DL architectures tuned and trained}
\label{app1}
All models were trained using Adam optimiser, with a batch size of 96, dropout rate of 0.1, and up to 500 epochs. The EarlyStopping method was implemented to prevent overfitting, and it allowed the training process to stop automatically if the validation loss does not improve over 50 consecutive epochs.
For CNN layers in Attention-based CNN-LSTM models, the kernel size was set to 3, pool size to 2, padding to `same' and activation function to Relu. The activation function of LSTM layers of the same model was set to Relu as well.
\begin{table}[ht]
    \centering
    \begin{tabularx}{0.48\textwidth}{|c|X|}
    \hline
        Model & Layers \\
    \hline
         \multirow{6}{*}{simpleRNN} &  Conv1D(filters=32, kernel size=24*4)   \\
         \cline{2-2}
         
                   & SimpleRNN(filters=64) \\
        \cline{2-2}
                  & SimpleRNN(filters=128) \\
        \cline{2-2}
                  & Flatten \\
        \cline{2-2}
                  & RepeatVector(96) \\
        \cline{2-2}
                  & TimeDistributed(Dense: filters=3) \\
    \hline
        \multirow{8}{*}{LSTM} &  Conv1D(filters=32, kernel size-24*4)   \\
         \cline{2-2}
                   & Conv1D(filters=63, kernel size=1) \\
        \cline{2-2}
                   & LSTM(filter=32) \\
        \cline{2-2}
                   & Flatten \\
        \cline{2-2}
                   &  RepeatVector(96)\\
        \cline{2-2}
                   & LSTM(filter=64) \\
        \cline{2-2}
                   &  LSTM(filter=32)\\
        \cline{2-2}
                   &  TimeDistributed(Dense: filter=3)\\
    \hline
         \multirow{9}{*}{bi-LSTM}  & Conv1D(filters=32, kernl size=24*4)  \\
         \cline{2-2}
         & Conv1D(filters=64, kernel size=1) \\
         \cline{2-2}
         & Bidirectional(LSTM: filters=128) \\
        \cline{2-2}
        &  Bidirectional(LSTM: filters=64)\\\cline{2-2}
        &  RepeatVector(96) \\
        \cline{2-2}
        &  LSTM(filters=128)\\
        \cline{2-2}
        &  LSTM(filters=64) \\
        \cline{2-2}
        &  LSTM(filters=32) \\
        \cline{2-2}
        &  TimeDistributed(Dense: filter=3) \\
    \hline
         \multirow{6}{*}{GRU} & Conv1D(filters=32, kernel size=24*4)   \\
         \cline{2-2}
        &  GRU(filters=64) \\
        \cline{2-2}
        &   GRU(filters=128) \\
        \cline{2-2}
        &  flatten \\
        \cline{2-2}
        & RepeatVector(96) \\
        \cline{2-2}
        & TimeDistributed(Dense: filters=3) \\
    \hline 
        \multirow{7}{*}{Attention-based CNN-LSTM}  &  LSTM(filters=50) \\
        \cline{2-2}
        & Attention layer \\
        \cline{2-2}
        &  RepeatVector(35) \\
        \cline{2-2}
        & Concatenate(LSTM, Attention) \\
        \cline{2-2}
        &  LSTM(filters=128) \\
        \cline{2-2}
        &  Dense(filters=96*3) \\
        \cline{2-2}
        &  Reshape(filters=96*3) \\
    \hline 
         Transformer & Transformer(layers=4, d\_model=128, heads=8, dff=512) \\
    \hline    
    \end{tabularx}
    \caption{DL architectures considered for the first round evaluation.}
    \label{tab:train_1}
\end{table}
%
\begin{table}[h]
    \centering
    \begin{tabular}{|c|c|l|}
    \hline
    Model & Parameter & Range/Set \\
    \hline
    Transformer & N. of layers & \{4, 5, 6\}  \\
    \cline{2-3}
                & d\_model &  \{32, 64, 128\} \\
    \cline{2-3}
                & N. of heads &  \{4, 8, 16\} \\
    \cline{2-3}
                & dff & \{128, 256, 512\}  \\
    \hline
    Attention-based CNN-LSTM & N. of CNN layers & \{1, 2\}  \\
    \cline{2-3}
                             & Filters & \{32, 64, 128\} \\
    \cline{2-3} 
                & N. of LSTM layers &  \{2, 3\} \\
    \cline{2-3}             
                 &  units & \{64, 128\} \\
    \hline
    \end{tabular}
    \caption{DL architectures and parameters tuned for the second round evaluation.}
    \label{tab:train_2}
\end{table}

\section*{Acknowledgements}
The authors would like to thank the funding for DriVe2X research and innovation project from the European Commission and the UKRI, with grant numbers 101056934 and 10055673, respectively.


\footnotesize
\balance
\bibliographystyle{IEEEtran} 
\bibliography{main}

\end{document}